\documentclass[11pt,twoside]{article}

\usepackage[utf8]{inputenc}
\usepackage{epsfig,amsfonts,color}
\usepackage{mathrsfs} 
\usepackage{amsmath, bbm, mathabx, amsthm}
\usepackage[normalem]{ulem}
\usepackage{amssymb, palatino, geometry,url}
\usepackage[colorlinks=true,linkcolor=blue,citecolor=blue,urlcolor=blue]{hyperref}
\usepackage{subcaption}
\usepackage{multirow}

\usepackage{graphicx} 

\usepackage{fancyhdr}
\usepackage{lineno}
\usepackage{float}
\usepackage[section]{placeins}
\usepackage[table]{xcolor}

\title{A Digital Twin Simulator of a Pastillation Process \\ with Applications to  Automatic Control\\ based on Computer Vision}

\author{Leonardo D. Gonz\'alez$^1$, Joshua L. Pulsipher$^2$, Shengli Jiang$^3$\\ Tyler Soderstrom$^4$, and Victor M. Zavala$^{1}$\thanks{Corresponding Author: zavalatejeda@wisc.edu}\\
    {\small $^1$Department of Chemical \& Biological Engineering}\\
    {\small \;University of Wisconsin-Madison, 1415 Engineering Drive, Madison, WI 53706, USA}\\
    {\small $^2$Department of Chemical Engineering}\\
    {\small \;University of Waterloo, 200 University Ave W, Waterloo, ON N2L 3G1, Canada}\\
    {\small $^3$Chemical \& Biological Engineering Department}\\
    {\small \;Princeton University, 50-70 Olden St, Princeton, NJ 08540, USA}\\
    {\small $^4$Advanced Process Control}\\
    {\small \;ExxonMobil Technology and Engineering, 22777 Springwoods Village Pkwy, Spring, TX 77389, USA}
    }
\date{}

\begin{document}

\maketitle

\begin{abstract}
We present a digital-twin simulator for a pastillation process. The simulation framework produces realistic thermal image data of the process that is used to train computer vision-based soft sensors based on convolutional neural networks (CNNs); the soft sensors produce output signals for temperature and product flow rate that enable real-time monitoring and feedback control. Pastillation technologies are high-throughput devices that are used in a broad range of industries; these processes face operational challenges such as real-time identification of clog locations (faults) in the rotating shell and the automatic, real-time adjustment of conveyor belt speed and operating conditions to stabilize output. The proposed simulator is able to capture this behavior and generates realistic data that can be used to benchmark different algorithms for image processing and different control architectures. We present a case study to illustrate the capabilities; the study explores behavior over a range of equipment sizes, clog locations, and clog duration. A feedback controller (tuned using Bayesian optimization) is used to adjust the conveyor belt speed based on the CNN output signal to achieve the desired process outputs.
\end{abstract}

\section{Introduction}

Pastillation processes are an important solidification technology that convert a molten product (i.e., polymer melts) into dispersed hemispherical solid forms called pastilles. These high-throughput processes are able to efficiently produce large quantities of uniformly-sized pastilles. Moreover, they offer a dust-free alternative to traditional mechanical cutting and breaking processes, which is particularly beneficial in environments where dust is problematic and/or when a hazardous chemical is used \cite{kim2003prediction, shukla2011lipid}. These benefits have led to the widespread adoption of pastillation in the chemical, pharmaceutical, food, and cosmetic industries \cite{shukla2011lipid, chen2017continuous}. For instance, in pharmaceutical manufacturing, direct solidification is crucial for producing drug products with precise control over droplet shape, necessary for enhancing drug release rates. Tablets made from crystalline materials can more effectively deliver drugs with controlled release rates compared to those made from amorphous materials. Additionally, because pastillation is a solvent-free process, direct solidification results in an exceptionally pure product \cite{guirgis2001hot, guirgis2001hot2}.
\vspace{0.1in}

\begin{figure}[!htb]
    \centering
    \includegraphics[width=\textwidth]{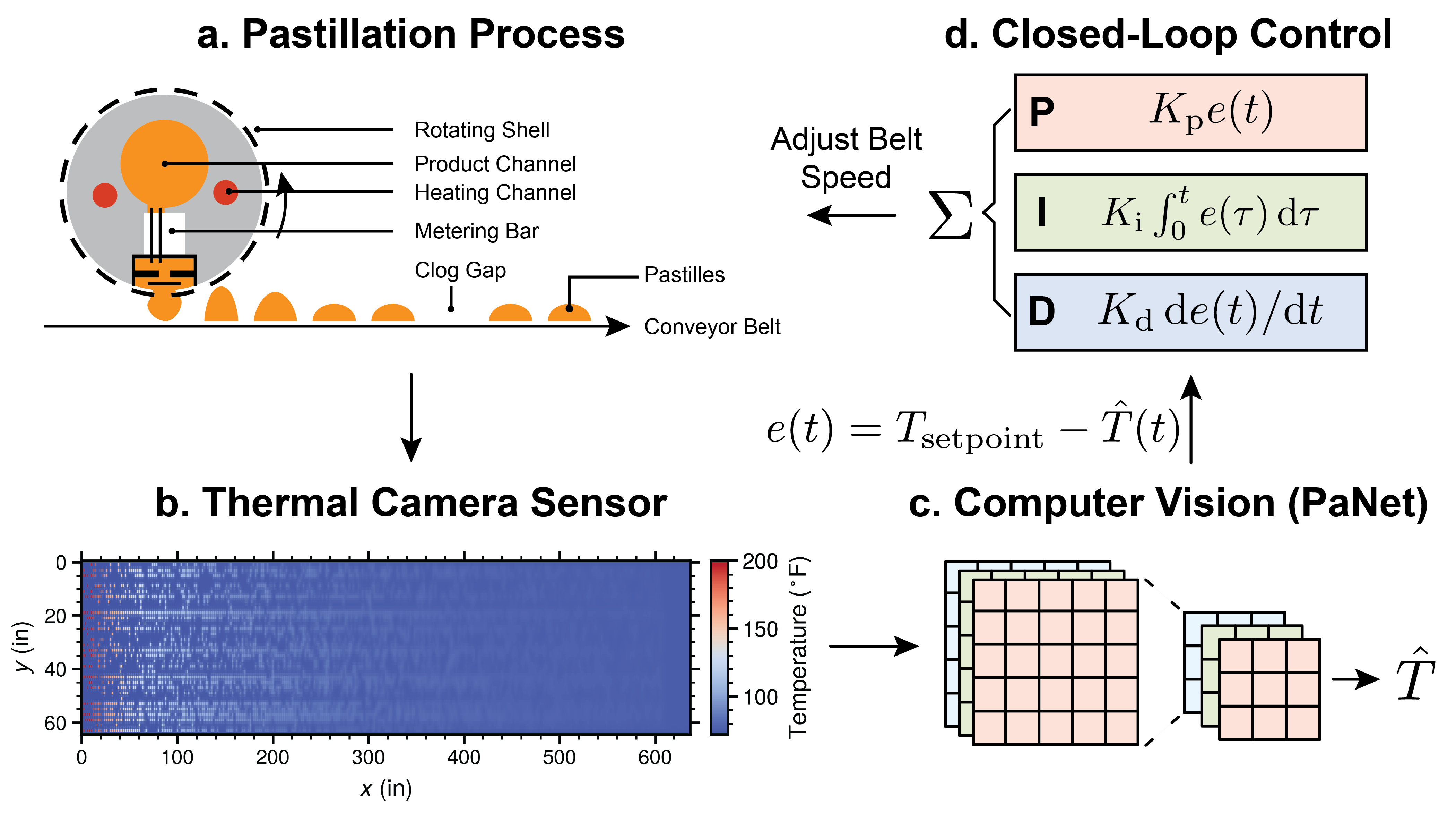}
    \caption{Overview of digital twin simulator for pastillation process.  
    (a) Schematic of a typical pastillation process: a rotating shell deposits pastilles onto a conveyor belt; 
    random clogging events introduce product flow variations.
    (b) A thermal camera sensor captures spatial temperature distributions of the conveyor belt in real-time; the simulator used models to generate realistic image data.  
    (c) A trained computer vision model (the convolutional neural network \texttt{PaNet}) processes thermal images to predict the average instantaneous flow and temperature of the pastilles.  
    (d) A feedback controller adjusts the conveyor speed to maintain the desired flow and temperature.}
    \label{fig:overview}
\end{figure}

A high-level schematic of a prototypical pastillation process is shown in Figure \ref{fig:overview}a. A melt product is fed into a heated cylindrical unit (e.g., a Rotoformer \cite{IPCO_2024}). Droplets of the melt are extruded through a rotating perforated shell on the unit, and these droplets are deposited on a cooled steel conveyor belt (which typically matches the rotation speed of the rotating shell) where they solidify into pastilles. The pastilles are then collected for packaging at the end of the conveyor belt. Physical properties of the melt that have direct implications for the size and shape of pastilles include viscosity, surface tension, density, and crystallization behavior. Process parameters include dropping height, heating/cooling rates, conveyor belt speed, and the product inlet pressure.
\\

A common operational challenge of pastillation processes is that nozzles on the rotating shell will randomly clog (fully or partially); this induces fluctuation in pastille production rates and sizing. To compensate for this fluctuation, operators typically adjust process parameters, namely the belt speed and inlet pressure, based on manual visual inspections of the conveyor belt (which provides qualitative measurement). However, the frequency and effectiveness of these manual interventions are limited due to competing operator tasks and the lack of real-time quantitative measurement of pastille production rates. Another key challenge pertains the temperature control of the pastilles, because visual cues are more subtle. The development of real-time sensing technologies that can provide information on flow rate and temperature of pastilles would facilitate monitoring of production and would be key enabler for automatic control to increase process consistency and efficiency.
\\

Drawing inspiration from how human operators currently assess pastille production rate via visual cues, automated computer vision provides an approach to measure production rates based on image data collected using high-speed cameras. In this context, computer vision refers to the use of data processing techniques (such as machine learning models) to make predictions based on image/video data \cite{voulodimos2018deep}. Convolutional neural networks (CNNs) are a popular choice of computer vision model that aim to automatically extract feature information from image data \cite{jiang2021convolutional, jiang2024convolutional, yoo2015deep}. In other words, CNNs provide a mechanism to map real-time image data to observable signals (thus generating a soft sensor); the observable signals can then be used for monitoring and feedback control \cite{gao2022augmented}. Computer vision sensors are widely used in robotics \cite{wiriyathammabhum2016computer}, autonomous vehicles \cite{janai2020computer}, consumer electronics \cite{alyamkin2019low}, healthcare \cite{gao2018computer}, manufacturing \cite{kakani2020critical, zhou2022computer}, and many other application areas. In manufacturing, computer vision sensors have been applied to conveyor belt units to measure diverse quantities such as belt speed \cite{gao2019contactless}, belt loading \cite{zhang2020computer}, belt health \cite{chamorro2022health}, and material flow rate \cite{lange1990measurement, sabih2023raw}. However, to the best of our knowledge, computer vision sensors have not been developed for pastillation processes.
\\

The utility of computer vision sensors in facilitating automatic control is well-established in traditional application areas such as autonomous vehicles \cite{bertozzi2000vision}. Moreover, such sensors are increasingly deployed in manufacturing to automate processes that have historically relied on operators making manual control actions based on visual information (e.g., a video feed) \cite{pulsipher2022safe}. These deployments reap the benefits of automatic control which include increased product consistency, decreased operator load, and improved process safety \cite{bequette2003process}. For instance, such systems have enabled effective automatic control systems for food drying \cite{martynenko2017computer}, printing cylinders \cite{villalba2019deep}, and chemical extrusion \cite{kadam2022systems}. {\em A practical limitation in developing computer vision systems for feedback control} is the lack of simulation environments (digital twins) that are capable of systematically generating data for training and validating machine learning and feedback control architectures. Such capabilities can help accelerate deployment of intelligent automation systems (e.g., via testing of technologies or development of transfer learning methodologies).
\\

In this work, we present a simulation framework of a pastillation process. The framework generates realistic image data of the conveyor belt that is fed to a computer vision  system to generate output signals that are in turn used to automatically adjust conveyor belt conditions using feedback control. Figure~\ref{fig:overview} provides an schematic of the simulation framework, which can be seen as a digital twin of the pastillation process (it captures sensors and control elements). The simulator models the spatial distribution, cooling, and production rate of pastilles on a conveyor belt subject to random clogging of the extrusion nozzles, and it produces top-view image data of the conveyor belt to simulate an overhead thermal camera. The simulator enables the generation of rich datasets that are used to train a computer vision sensor (that we call \texttt{PaNet}-Pastillation Network); this soft sensor uses a CNN architecture to predict the instantaneous average flow rate and temperature of the pastilles. Finally, we develop a closed-loop control system that uses \texttt{PaNet} and proportional-integral-derivative (PID) controllers to control the temperature and flow rate of pastilles. We demonstrate the effectiveness of this approach relative to manual control via simulation studies.
\\

The paper is structured as follows. Section \ref{sec:simulator} describes the proposed pastillation simulator. Section \ref{sec:sensor} details the proposed CNN sensor for computer vision. Section \ref{sec:control} defines the proposed feedback control system. Section \ref{sec:results} presents and discusses simulation results, and Section \ref{sec:conclusions} provides concluding remarks and directions for future work.

\begin{figure}[htbp]
    \centering
    \includegraphics[width=0.80\textwidth]{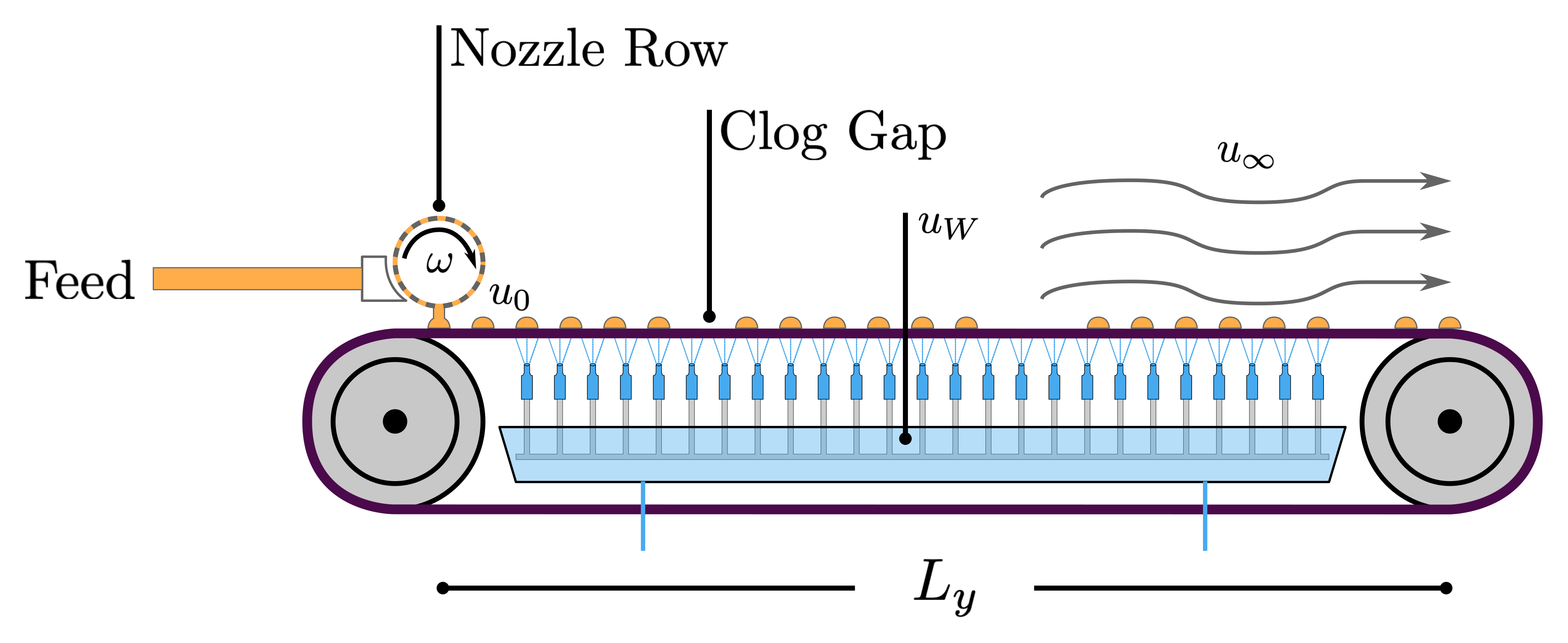}
    \caption{Schematic representation of pastillation system. A rotating shell with speed $\omega$ and K evenly-spaced rows of $h$ nozzles extrudes pastilles at a temperature $u_0$ along the width, $L_x$ of a belt. As they move down the length, $L_y$, of the belt, the pastilles are cooled via cooling water sprayed onto the underside of the belt (at temperature $u_W$) and the ambient air (at temperature $u_{\infty}$).}
    \label{fig:pastillation_belt}
\end{figure}

\section{Pastillation Process Simulator} \label{sec:simulator}

We detail the development of the pastillation conveyor belt simulator, which models both the flowrate and the temperature field of the pastilles as they move down the conveyor belt. The pastillation system, shown in Figure \ref{fig:pastillation_belt}, consists of a rotating shell moving at an angular speed $\omega$ with $H$ nozzles configured in $K$ rows with $h=\frac{H}{K}$ nozzles per row. A new row of $h$ pastilles is placed onto the belt with a frequency of $\frac{\omega K}{2\pi}$. 
\\

The extruded droplets are evenly spaced along the width of the belt ($L_x$) and exit the nozzles at temperature $u_0$. As the pastilles travel down the length of the belt ($L_y$), they are cooled by convection with ambient air (at temperature $u_\infty$) and by a system of water jets that spray cooling water (at temperature $u_{W}$) onto the underside of the belt. 
\\

A detailed list of the values of relevant model parameters along with the computer code used to construct the simulator can be found at  \url{https://github.com/zavalab/ML/tree/master/Pastillation}. 

\begin{figure}[htbp]
    \centering
    \includegraphics[width=0.7\textwidth]{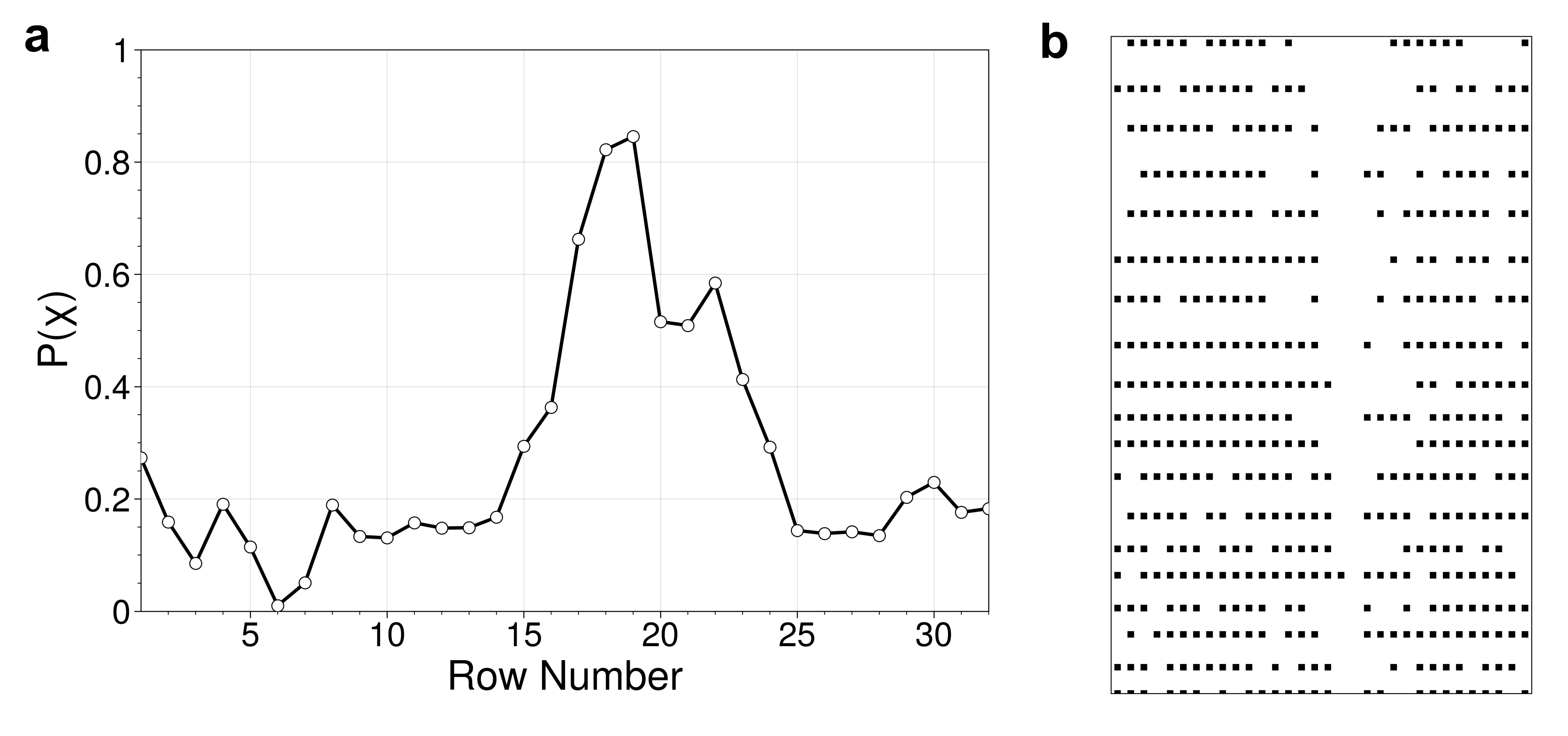}
    \caption{Simulation of nozzle clogs. (a) Empirical probability distribution used for determining the occurrence of a clog at each nozzle. (b) Representative distribution of pastille placement onto the belt.}
    \label{fig:clog_and_flow_sample}
\end{figure}

\subsection{Modeling of Pastille Flow Rate}

Maintaining a consistent pastille flow rate (product flow rate) in an actual process is challenging due to the frequent and random formation of clogs in the extrusion outlets. These blockages are caused by a wide array of factors such as extrusion speed, temperature gradient, and melt thickness, which are difficult to model and predict. Moreover, 
some rotoformer nozzles can exhibit a higher tendency to block than others, making it inaccurate to assign a uniform clog probability value. 
\\

We used real image data from an industrial system to generate an empirical clog probability distribution based on the clog gaps observed for each nozzle across various pastille rows. The distribution, shown in Figure \ref{fig:clog_and_flow_sample}a, was constructed using a scaled temperature where the probability of a clog occurring, $P(\chi)$, was defined as 
\begin{equation}
    P(\chi)=1-\frac{u-u_\infty}{u_{M}-u_\infty}
\end{equation}
where $u_{M}$ is the maximum temperature measured in the images and $u$ is the average temperature observed for a row at the head of the belt. This approach allows for an accurate simulation of the clogging patterns of the extrusion system and the resulting random fluctuations in the production rate, as seen in Figure \ref{fig:clog_and_flow_sample}b. Whenever the rotating shell places a new row of pastilles onto the belt, the generated distribution is sampled to determine if a clog occurs at any of the outlets (on such row). The actual number of pastilles placed is then simply the difference between $h$ and the number of blocked nozzles. Each row is tracked as it moves down the belt in order to calculate the total number of pastilles on the belt at any given time ($\Theta$). Given a belt speed $v_{B}$, the average production rate of the system $\dot{m}_{P}$ can be determined as follows:

\begin{equation}
    \dot{m}_{P} = \frac{\Theta}{L_y}\cdot v_{B}
\end{equation}

\subsection{Modeling of Spatial Temperature Field}

The dynamic evolution of the spatial temperature field of the conveyor belts is modeled using a 2D heat equation:
\begin{subequations}\label{eq:2d_heat}
\begin{align}
    & \frac{\partial u}{\partial t}=\alpha\left(\frac{\partial^2u}{\partial x^2}+\frac{\partial^2u}{\partial y^2}\right)+f(t, x, y) \label{eq:heat_pde} \\
    & u(t, x=0, y) = u_{\infty}, ~~~ u(t, x=L_x, y) = u_{\infty}\\
    & u(t, x, y=0) = u_{\infty}, ~~~ u(t, x, y=L_y) = u_{\infty}\\
    & u(t=0, x, y) = I(x, y)
\end{align}
\end{subequations}
Here, $u(t, x, y)$ is the temperature at time $t$ and location $(x, y)$ on the belt; $\alpha$ is the thermal diffusivity, which measures how quickly heat moves through the material; and $f(t, x, y)$ is a forcing function that captures the effects of any heat sources/sinks present within the system. 
\\

The above model is defined over the domain $[0, T]\times[0, L_x]\times[0, L_y]$, where $T$ is the time duration of the simulation. The boundary conditions on all sides are given by $u_{\infty}$, and the initial temperature distribution of the system is described by $I(x,y)$. 
\\

The thermal model is solved using a finite difference scheme \cite{langtangen2017finite}. This approach discretizes the spatiotemporal domain into a mesh of $N_t\times N_x\times N_y$ elements. Using finite difference schemes, \eqref{eq:2d_heat} is converted into the following set of algebraic equations:
\begin{subequations}\label{eq:finite_difference}
\begin{gather}
    \frac{u_{i,j}^{n+1}-u_{i,j}^{n}}{\Delta t}=\alpha\left(\frac{u_{i-1, j}^{n+1}-2u_{i, j}^{n+1}+u_{i+1, j}^{n+1}}{\Delta x^2}+\frac{u_{i, j-1}^{n+1}-2u_{i, j}^{n+1}+u_{i, j+1}^{n+1}}{\Delta y^2}\right)+f^n_{i,j}\\
    \Delta t=\frac{T}{N_t},~~\Delta x=\frac{L_x}{N_x},~~\Delta y=\frac{L_y}{N_y}
\end{gather}
\end{subequations}
which is defined over the domain $i\in[1,N_x-1]$, $j\in[1,N_y-1]$, and $n\in[0,N_t-1]$. Grouping like terms then yields:
\begin{subequations}\label{eq:algebraic_system}
\begin{gather}
    (1+2(F_x+F_y))u_{i,j}^{n+1}-F_x(u_{i-1,j}^{n+1}+u_{i+1,j}^{n+1})-F_y(u_{i,j-1}^{n+1}+u_{i,j+1}^{n+1})=u_{i,j}^{n}+f_{i,j}^{n}\Delta t\\
    F_x=\frac{\alpha\Delta t}{\Delta x^2},~~F_y=\frac{\alpha\Delta t}{\Delta y^2}
\end{gather}
\end{subequations}
where $F_x$ and $F_y$ are referred to as the Fourier numbers in the $x$ and $y$ directions, respectively. At each time step $n$, a system of linear equations of the form $Aw=b$ is constructed by aggregating the corresponding algebraic representations of \eqref{eq:2d_heat} at each mesh point in $[1,N_x-1]\times[1,N_y-1]$. 
\\

The temperature within each element at the next time step, $n+1$, is then determined by solving this system for $w$. Once the solution has been obtained, the $y$ positions of each row are updated based on the belt speed as follows:
\begin{equation}
    j^{n+1}=j^n+\frac{v_B\Delta t}{\Delta y}
\end{equation}
Note that once $j>N_y$ for any row, this has moved past the end of the belt and is no longer tracked. This process is repeated recursively across the specified simulation time $T$ to obtain the dynamic profile of the temperature field.
\\

The effect of the cooling water system on the temperature field is captured via the forcing function $f(t, x, y)$. This system is comprised of a set of water jets uniformly distributed along $R$ rows with $q$ jets per row, each spraying cooling water onto the underside of the belt at a rate $\dot{m}_{W}$. Each row wets an area $S_r, ~ r=1,...,R$, centered around the jets, with every point within this region receiving the same water flow rate (i.e., the spray is uniformly dispersed). While the belt has a thickness $\ell$, this is small enough so that the thermal resistance between the top and bottoms sides of the belt is negligible. It is also assumed that the contact time between the cooling water and the belt is sufficient to reach thermal equilibrium. Thus, following the form of \eqref{eq:algebraic_system}, the forcing function is defined as:
\begin{equation}\label{eq:forcing_function}
    f^n_{i,j}=
        \begin{cases}
            -\frac{q\dot{m}_{W}C_{p_{W}}}{S_R\ell\rho_BC_{p_B}}(u^{n}_{i,j}-u_{W}) & \textrm{if} ~~ (i,j) \in S_r, ~ r = 1,...,R\\
            0 & \textrm{otherwise}
        \end{cases}
\end{equation}
Here $C_{p_{W}}$ and $C_{p_{B}}$ are the heat capacities of the cooling water and belt, respectively, and $\rho_B$ is the density of the belt. It should be noted that these physical parameters, as well as $\alpha$, are assumed to be invariant with time or position. However, the modeling framework can be easily extended/refined to account for more sophisticated behavior.

\section{Computer Vision Sensor (\texttt{PaNet})} \label{sec:sensor}

In this section, we introduce \texttt{PaNet}, a computer vision sensor that is designed to predict average temperature and flow rate from thermal image data. We outline the data preprocessing steps, model architecture, and training strategies used to optimize its performance.

\subsection{Data Preprocessing}\label{sec:data_prep}

To improve the generalizability of \texttt{PaNet} for real-world pastillation processes, we curated a diverse dataset obtained from simulated snapshots. Each snapshot is a single-channel image, simulating a thermal camera, with pixel values normalized between 0 and 1. This scale corresponds to a temperature range of $72^{\circ} \mathrm{F}$ (background) to $212^{\circ} \mathrm{F}$, adjustable to specific pastillation processes. This range was used as an example, and its adjustment to other values does not impact algorithm performance. The images are represented as $\mathbf{X} \in \mathbb{R}^{637 \times 65}$, where 637 pixels represent the length of the conveyor belt and 65 pixels represent its width. The corresponding label, $\mathbf{y} \in \mathbb{R}^2$, encodes both the average temperature of the earliest produced pastilles in the image and the pastille flow rate ($\dot{m}_P$). Early analysis of pastilles, even before they traverse the full conveyor belt, is critical for detecting process anomalies such as slow initial speeds that cause overcooling. Using the entire image rather than just the final segment enables early detection of such issues, allowing for rapid process adjustments.

\subsection{Model Architecture}\label{sec:mod_arc}

To predict both the average temperature and the pastille flow rate, we developed a couple of CNN architectures: a 2D CNN (\texttt{PaNet-2D}) and a 1D version (\texttt{PaNet-1D}). Both models take a single-channel image of size $637 \times 65$ as input, where 637 represents the length of the conveyor belt and 65 represents its width. In \texttt{PaNet-2D}, 2D convolutional layers with $3 \times 3$ filters extract patterns from the input. Each convolution is followed by a $2 \times 2$ average-pooling operation to reduce spatial dimensions. The architecture includes three convolutional layers, three average-pooling layers, and three fully connected layers. In \texttt{PaNet-1D}, the input image is treated as 65 channels, where each channel corresponds to a width vector of size 637. The architecture consists of three 1D convolutional layers with $3$ filters, 3 1D average-pooling layers, and 3 fully connected layers. In both models, the hyperparameter $n$ determines the number of filters in convolutional layers and hidden units in fully connected layers. A ReLU activation is applied throughout, including the output layer, to ensure non-negative predictions.
\\

\subsection{Model Training and Hyperparameter Tuning}\label{sec:model_train_hyper}

To evaluate model generalizability and minimize selection bias, we employed five-fold cross-validation (CV) on a dataset that contained 20,000 images. The data was split into training, validation, and test sets containing 12,800, 3,200, and 4,000 images, respectively. For 5-fold CV, the dataset was divided into five equal subsets. In each iteration, one subset was used as the test set, and the remaining four served as the training set. This process was repeated five times, ensuring each subset was used as the test set once.
\\

The training set was used to fit the model, with 20\% reserved as a validation set for hyperparameter tuning and overfitting prevention. Model performance was evaluated on the independent test set. Hyperparameter tuning included batch sizes $\{32, 64, 128\}$, learning rates $\{0.0005, 0.001, 0.005\}$, and filter and hidden unit counts $\{64, 128, 256\}$. The loss function was mean squared error (MSE), and optimal hyperparameters for each fold were selected based on the lowest average validation root mean squared error (RMSE). This approach yielded the mean and standard deviation of RMSE and $R^2$ scores across all 5 test sets, ensuring robust evaluation of model accuracy.

\subsection{Saliency Analysis}\label{sec:saliency_analysis}

We employed the Python package {\tt SmoothGrad}~\cite{smilkov2017smoothgrad} to obtain saliency maps for interpreting the model predictions. Let $f(\mathbf{x})$ denote our CNN model, which generates a prediction $\hat{\mathbf{y}}\in\mathbb{R}^2$ for an input image $\mathbf{x}$, corresponding to the temperature and flow rate. During training, we minimize the mean squared error (MSE) loss,
\begin{equation}\label{eq:mse}
\mathcal{L}(\mathbf{y}, \hat{\mathbf{y}}) = \frac{1}{N} \sum_{i=1}^{N} (\mathbf{y}_i - \hat{\mathbf{y}}_i)^2,
\end{equation}
where $\mathbf{y}_i$ contains both the true temperature and flow rate, and $\hat{\mathbf{y}}_i$ represents the corresponding predicted values for the $i$th sample. For a single input $\mathbf{x}$, we define the saliency map as the gradient of the loss with respect to the input:
\begin{equation}\label{eq:gradient}
\mathbf{G}(\mathbf{x}) = \nabla_{\mathbf{x}} \mathcal{L}\bigl(\mathbf{y}, f(\mathbf{x})\bigr),
\end{equation}
which has the same dimensions as $\mathbf{x}$. In our setup, $\mathbf{x} \in \mathbb{R}^{637 \times 65}$ for \texttt{PaNet-1D} and $\mathbf{x} \in \mathbb{R}^{637 \times 65 \times 1}$ for \texttt{PaNet-2D}.
\\

Simple gradient-based saliency maps can be noisy. To address this, we apply {\tt SmoothGrad}, which averages gradients over multiple noisy perturbations of $\mathbf{x}$:
\begin{equation}\label{eq:smoothgrad}
\text{Saliency}_{\text{SmoothGrad}}(\mathbf{x})
= \frac{1}{M} \sum_{m=1}^{M} \nabla_{\mathbf{x}} \mathcal{L}\!\Bigl(\mathbf{y}, f\bigl(\mathbf{x} + \boldsymbol{\epsilon}_m\bigr)\Bigr),
\end{equation}
where $\boldsymbol{\epsilon}_m \sim \mathcal{N}(\mathbf{0}, \sigma^2 \mathbf{I})$ is noise sampled from a Gaussian distribution with mean 0 and standard deviation $\sigma = 0.001$. By introducing noise, {\tt SmoothGrad} reduces high-frequency artifacts, making the resulting saliency maps more robust and interpretable.
\\

After computing the {\tt SmoothGrad} map, we take its absolute value and normalize it via min-max scaling:
\begin{equation}\label{eq:normalize}
\widetilde{\mathbf{S}}(i,j) = \frac{\bigl|\text{Saliency}_{\text{SmoothGrad}}(\mathbf{x})(i,j)\bigr| - \min(\mathbf{S})}{\max(\mathbf{S}) - \min(\mathbf{S})},
\end{equation}
where $\mathbf{S} = \bigl\{\bigl|\text{Saliency}_{\text{SmoothGrad}}(\mathbf{x})(i,j)\bigr|\bigr\}$ over all pixels $(i,j)$ in the saliency map. This produces final saliency values in $[0,1]$, allowing for clear visualization of the most influential regions for temperature prediction.

\section{Feedback Control} \label{sec:control}

A PID controller was designed to regulate the temperature of the system by modulating the belt speed. As the controller requires a point measurement, we use the \texttt{PaNet} sensor to predict the  average temperature of the row furthest along the belt, which we also refer to as the {\em leading row}. This output value is passed onto the controller and used to determine the next control action. We note that the \texttt{PaNet} updates the position of the leading row at every time step. For any time step $n$, once the current temperature measurement, $U_n$, has been determined by \texttt{PaNet}, the control error is calculated as $e_n=\hat{U}-U_n$, where $\hat{U}$ is the temperature setpoint. This value is then appended to the error history, which is used to determine the adjustment of the belt speed, $s_n$. We used a controller of the form:
\begin{equation}
    \Delta s_n = K_P\left(e_n-e_{n-1}+\frac{\Delta t}{\tau_I}e_n+\tau_D\left(\frac{e_n-2e_{n-1}+e_{n-2}}{\Delta t}\right)\right)
\end{equation}
Here, $\Delta s_n$ is the change in belt speed at time $n$ and $\Delta t$ is the time step size. The controller gain is denoted by $K_P$ and $\tau_I$ and $\tau_D$ are the integral and derivative time constants, respectively. The subsequent speed value is then simply
\begin{equation}\label{eq:new_action}
    s_n = s_{n-1}+\Delta s_n.
\end{equation}
To prevent the controller from applying overly aggressive corrections (which can result in instability and cause significant wear on the system actuators), the change in the belt speed was bounded. This provides the system with additional time to respond more smoothly to the applied control action. The speed of the belt itself was also constrained such that $s_n\in[2, 12]$. Thus, updating \eqref{eq:new_action} based on these considerations results in the speed being set as follows:
\begin{subequations}
\begin{align}
    \Delta s_n &= \min(\max(\Delta s_n, -1), 1)\\
    s_n &= \min(\max(s_{n-1}+\Delta s_n,2),12).
\end{align}
\end{subequations}
The selection of the controller parameters $K_P$, $\tau_I$, and $\tau_D$ has a significant impact on the performance of the controllers. In practice, these are usually tuned using heuristics or expert knowledge. However, these methods can require a significant amount of trial-and-error and can result in suboptimal performance due to complex interactions between the controller parameters and due to the presence of constraints. To accelerate the tuning process, we automated the controller tuning using the Variable Partitioning Bayesian Optimization (VP-BO) algorithm \cite{gonzalez2023vpbo}. For each set of controller parameters proposed by the algorithm, the simulation was executed for $T=400$ time steps, and the sum of the absolute value of the set-point error at each time step was used as the performance metric. The goal of the algorithm was to solve the following problem:
\begin{subequations}\label{eq:opt_prob}
\begin{gather}
    \min_x~~\sum_{n=1}^{T}|e_n(x)|\label{eq:goal}\\
    s.t.~~x\in X
\end{gather}
\end{subequations}
where $x=[K_P, \tau_I, \tau_D]^T$ and $X$ is the domain or search space from which the solution is sought. VP-BO utilizes the generated input-output data to construct a surrogate model that estimates the relation between the controller parameters and the closed-loop performance metric. In addition to predicting performance values $\mathcal{E}(x)$ this model also measures the uncertainty of these predictions $\mathcal{U}(x)$. These metrics are combined in a utility or acquisition function:
\begin{equation}
\mathcal{AF}=\alpha(\mathcal{E}(x), \mathcal{U}(x);\kappa)
\end{equation}
where $\kappa$ is the exploratory weight that scales the value placed on the model uncertainty. A new set of parameters was then selected by solving the following problem
\begin{subequations}\label{eq:af}
\begin{gather}
    \min_x~~\mathcal{AF}(x;\kappa)\\
    s.t.~~x\in X
\end{gather}
\end{subequations}

This approach enables the algorithm to balance the value of sampling from unexplored regions of $X$ with refining observed promising regions. The controller parameters given by the solution of \eqref{eq:af} were tested in simulation and the resulting data was used to update the surrogate model. This process was repeated until the experiment budget of $B$ iterations was exhausted. The controller parameters that were observed to provide the best performance were selected.

\section{Case Study} \label{sec:results}

We demonstrate the capabilities of our simulation framework using a case study; the nature of the study was inspired by a real setting for an industrial pastillation process. 

\begin{figure}[!htbp]
    \centering
    \includegraphics[width=\textwidth]{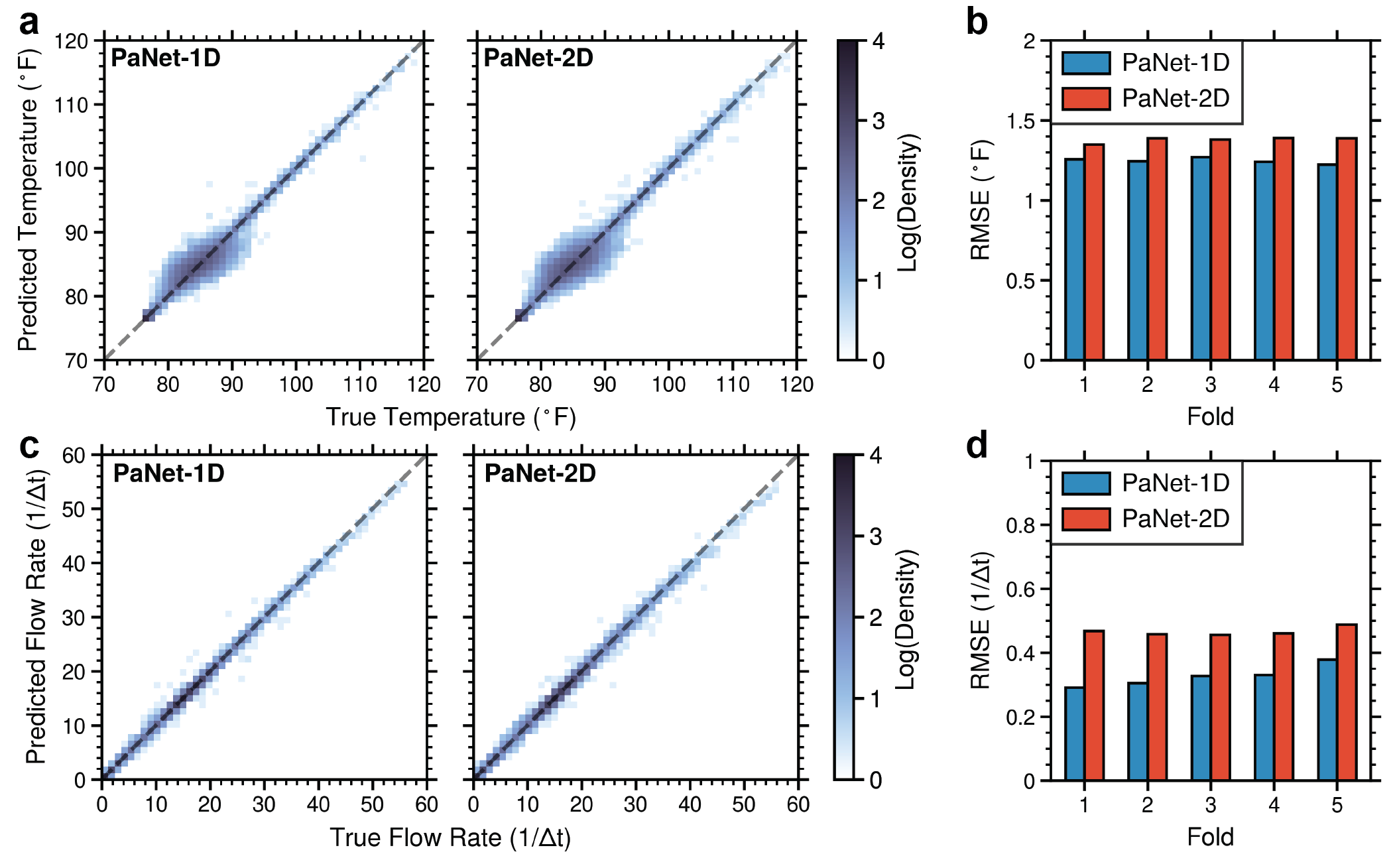}
    \caption{{\bf Performance of \texttt{PaNets} in predicting average temperature and pastille flow rate.} (a, c) Regression parity plots for \texttt{PaNet-1D} and \texttt{PaNet-2D}, with the diagonal line indicating ideal predictions. The colorbar represents the logarithm of density. (b, d) Root mean square error (RMSE) for each fold in five-fold cross-validation.}
    \label{fig:ml_result}
\end{figure}

\subsection{Computer Vision Sensor Performance} \label{sec:cv_result}

Figure~\ref{fig:ml_result} illustrates the performance of \texttt{PaNet-1D} and \texttt{PaNet-2D} in predicting average temperature and pastille flow rate. In Figures~\ref{fig:ml_result}a and c, the regression parity plots demonstrate that the CNN predictions closely match the observed outputs for both target quantities. Moreover, Figure~\ref{fig:ml_result}b reveal that both models achieve an average RMSE below $1.4^{\circ} \mathrm{F}$ across five folds of validation for temperature prediction. The \texttt{PaNet-1D} architecture outperforms \texttt{PaNet-2D}, with a mean RMSE of $1.2472 \pm 0.0157^{\circ} \mathrm{F}$ compared to $1.3794 \pm 0.0157^{\circ} \mathrm{F}$. Corresponding $R^2$ values are $0.9611 \pm 0.0020$ for the 1D model and $0.9525 \pm 0.0014$ for the 2D model. Similarly, for pastille flow rate prediction (Figure~\ref{fig:ml_result}d), \texttt{PaNet-1D} again outperforms \texttt{PaNet-2D}, achieving a mean RMSE of $0.3269 \pm 0.0297 /\Delta t $ compared to $0.4663 \pm 0.0117 /\Delta t$, with $R^2$ values of $0.9981 \pm 0.0003$ and $0.9961 \pm 0.0002$, respectively. These results indicate that the simpler 1D architecture effectively captures both thermal patterns and pastille distribution from pastillation images, providing a robust basis for feedback control of the pastillation processes.
\\

\begin{figure}[!htb]
    \centering
    \includegraphics[width=.85\textwidth]{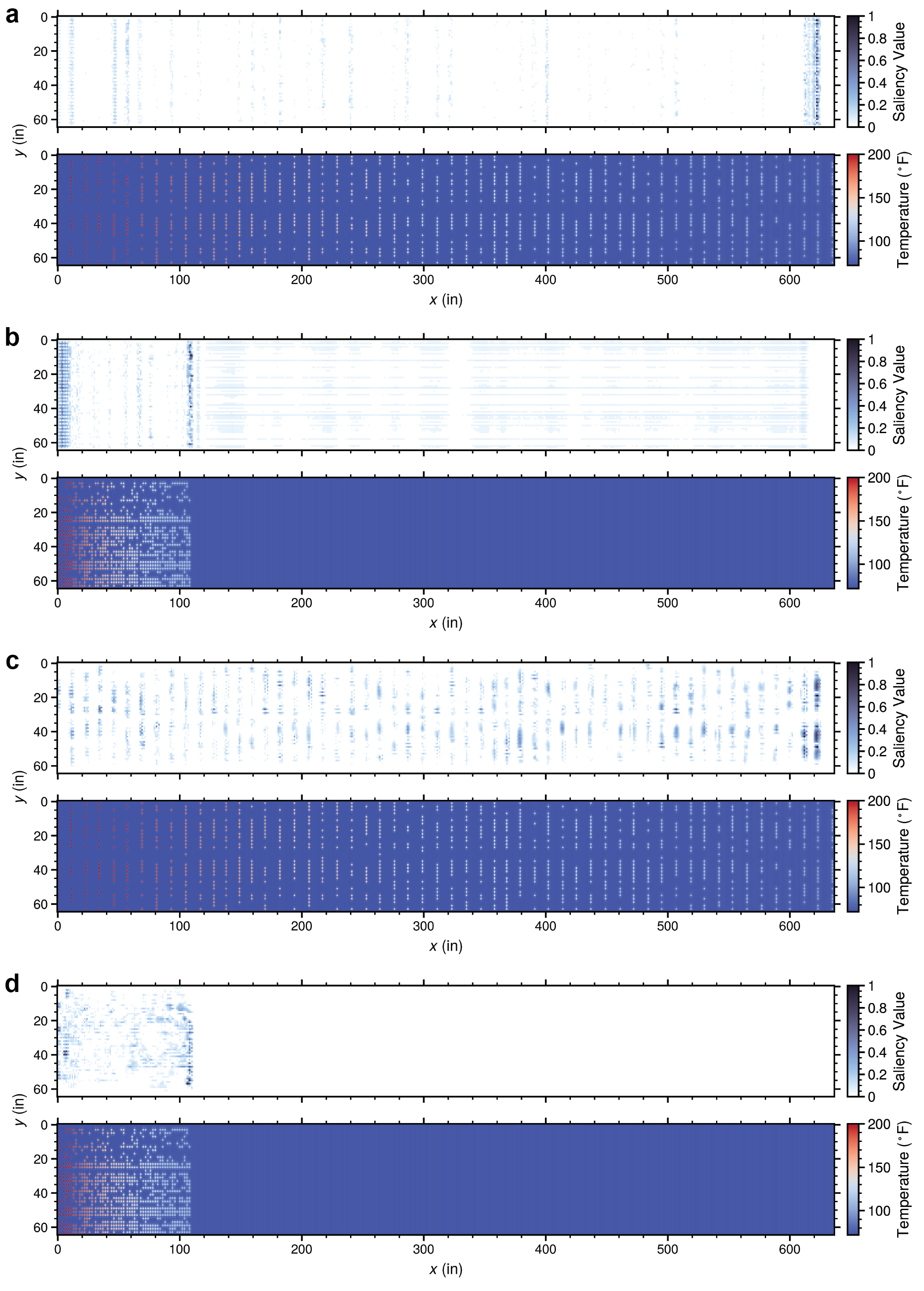}
    \caption{Saliency maps of \texttt{PaNets} for predicting average temperature. Panels (a) and (b) correspond to \texttt{PaNet-1D}, whereas (c) and (d) correspond to \texttt{PaNet-2D}. In (a) and (c), the earliest row of pastilles has reached the end of the conveyor belt; in (b) and (d), the earliest row is still in transit. For each panel, the top image shows the saliency map (darker regions indicate areas of greater relevance to the model), and the bottom image displays the corresponding thermal camera snapshot.}
    \label{fig:saliency}
\end{figure}

To gain additional insight into the prediction performance of \texttt{PaNet-1D}, we employed {\tt SmoothGrad} saliency maps to analyze how each model focuses on input data. In Figure~\ref{fig:saliency}a (\texttt{PaNet-1D}) and Figure~\ref{fig:saliency}c (\texttt{PaNet-2D}), the earliest row of pastilles has reached the end of the conveyor belt. Notably, \texttt{PaNet-1D} precisely highlights the final row in the thermal image, the region where temperature control is targeted, whereas \texttt{PaNet-2D} exhibits a more diffuse focus across the belt. In Figure~\ref{fig:saliency}b (\texttt{PaNet-1D}) and Figure~\ref{fig:saliency}d (\texttt{PaNet-2D}), the earliest row of pastilles remains in transit. Similarly, \texttt{PaNet-1D} accurately pinpoints the earliest row ($\sim 110$ inch), while \texttt{PaNet-2D} identifies the same row with less distinct focus. The targeted attention mechanism of \texttt{PaNet-1D} underpins its superior ability to track thermal patterns, enabling more accurate temperature predictions. This enhanced focus synergistically contributes to its high performance in simultaneous flow rate prediction. We note that saliency maps provide interesting information that could be used in the future to reduce the complexity of the CNN model; this can be done by identifying specific regions that have the largest impact on predictive performance. 

\vspace{0.1in}

\subsection{Feedback Control Performance} \label{sec:control_result}

We proceeded to connect the CNN sensor to the feedback controller. 
\\

\begin{figure}[!htb]
    \centering
    \includegraphics[width=\textwidth]{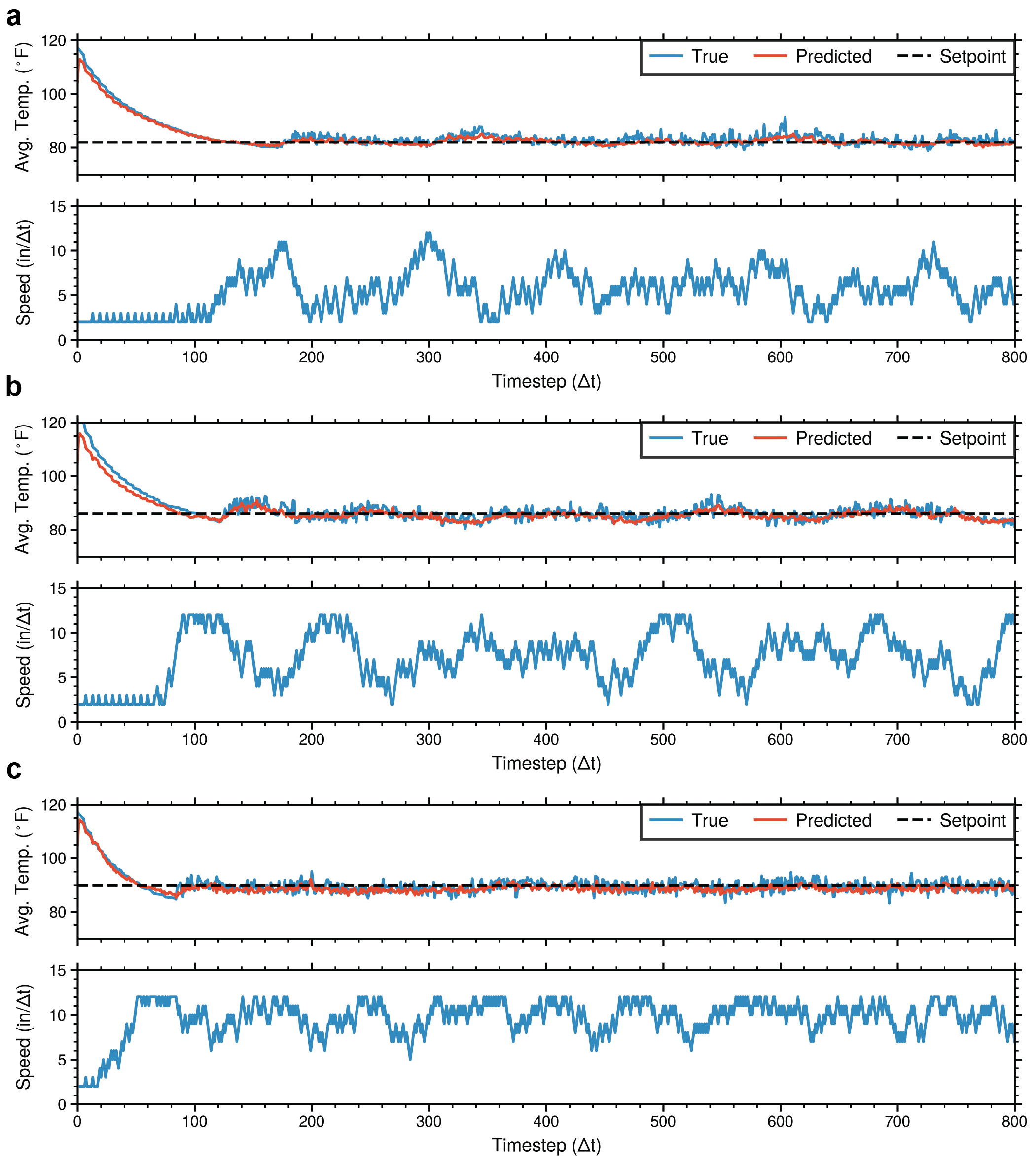}
    \caption{PID control system performance using \texttt{PaNet-1D} for temperature prediction. Panels (a), (b), and (c) show temperature control at setpoints of 82$^\circ$F, 86$^\circ$F, and 90$^\circ$F, respectively. In each panel, the top plot compares the true and predicted average temperatures of the earliest row of pastilles on the conveyor belt, with the black dashed line indicating the setpoint. The bottom plot shows the conveyor speed profile serving as the control variable.}
    \label{fig:control}
\end{figure}

Figure~\ref{fig:control} illustrates the control system performance using \texttt{PaNet-1D} for temperature prediction under three setpoint conditions (82$^\circ$F, 86$^\circ$F, and 90$^\circ$F) and using the optimal controller parameters provided by VP-BO ($K_P=47.0, \tau_I=15.3, \tau_D=0.0234$). Note that the setpoint temperature can be adjusted over wider ranges depending on the thermal camera scale and operating conditions. Overall, the predicted temperature of the earliest row of pastilles on the conveyor belt closely matches the true temperature for all setpoints. However, at the 90$^\circ$F setpoint, a larger discrepancy of about 5$^\circ$F occurs at high initial temperatures, possibly because of fewer representative high-temperature thermal images in the training set.
\\

All setpoints are reached within 150 timesteps ($\Delta t$), indicating the system responsiveness and efficiency. Notably, for the 90$^\circ$F setpoint, the target temperature is attained in fewer than 100 $\Delta t$. Meanwhile, the conveyor speed profile begins adjusting around 40 $\Delta t$ while the pastilles are still in transit, highlighting the system ability to anticipate and respond to early temperature deviations.
\\

Despite successful temperature control at all setpoints, both the true and predicted temperatures exhibit fluctuations around their targets. While the predicted temperature traces tend to be smoother, the true temperature readings show sharper, more rapid oscillations. Inspection of the belt speed, which is capped at 12 units per $\Delta t$, further reveals notable variations in the speed control profiles, particularly after convergence to the setpoint. These fluctuations suggest possible instability in the control approach, arising from model inaccuracies or overcompensation by the controller. Future refinements, such as advanced control algorithms, may mitigate these instabilities and further enhance system performance.
\vspace{0.1in}

\begin{figure}[htbp]
    \centering
    \includegraphics[width=1.0\textwidth]{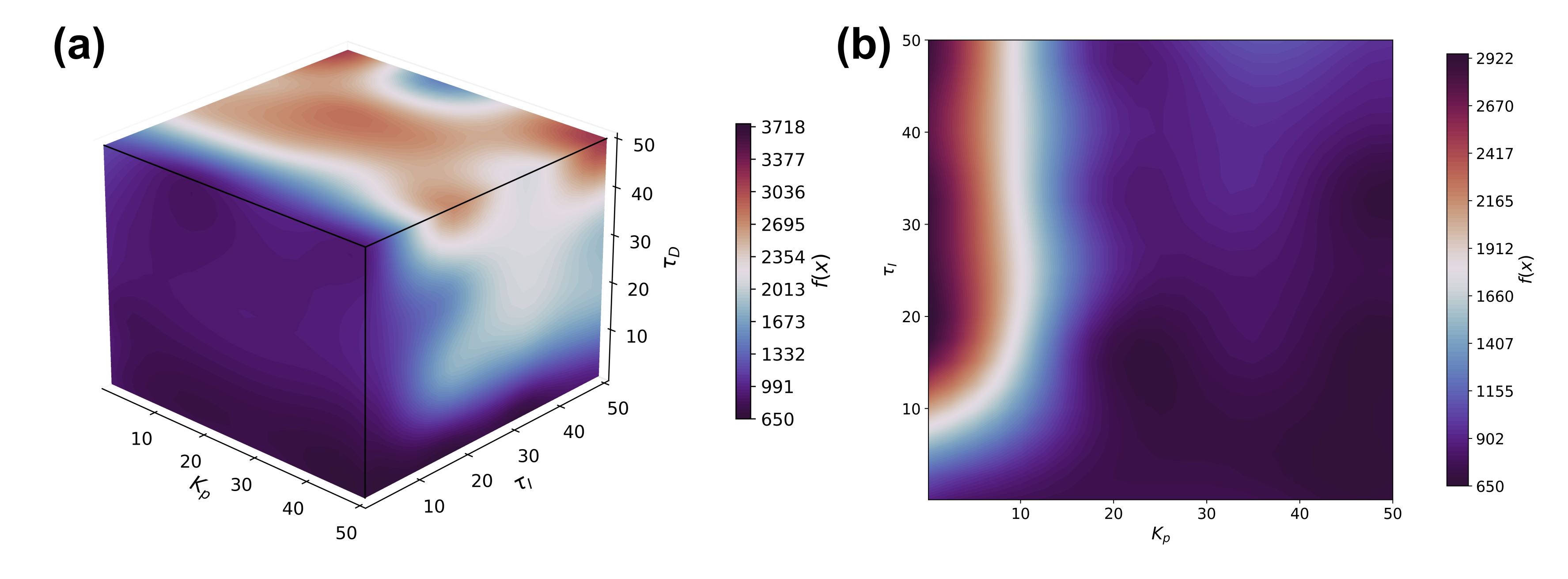}
    \caption{Controller tuning objective. (a) Box surface plot of tuning objective $f(x)$ in the defined control design space $X = [0.1,0.1,0.01]\times[50,50,50]$. (b) 2D slice of $f(x)$ at the optimal value of $\tau_D = 0.0234$.}
    \label{fig:ccontroller_tuning_objective}
\end{figure}

In order to illustrate the performance of the BO-based tuning algorithm, in Figure \ref{fig:ccontroller_tuning_objective} we show the controller performance as a function of the controller parameters. Here, we can clearly see that there is highly complex dependency of performance with respect to the parameters. The VP-BO tuning algorithm ran for 10 iterations across 3 partitions for a total of 30 experiments. A Gaussian process constructed from a M\'atern kernel with a smoothness parameter of 1.5 was used to generate the surrogate model. The lower confidence bound (LCB) function was used as the acquisition function, and the exploratory parameter was set to $\kappa = 2.6$. All simulations were conducted using a temperature setpoint of $\hat{U} = 90$. Additionally, the starting index in \eqref{eq:goal} was shifted from 1 to the time step at which the first row placed on the belt exits the system, $t_{N_y}$. This was due to the significant error that is observed during the system initialization. This is a purely transient event that is a direct result of the pastilles exiting the rotating shell at a very high temperature and is not representative of long-term system behavior. Thus, to avoid the potential biasing effects of the large error observed during the initial simulation period, these values were not considered in the calculation of the performance.  

\section{Conclusions} \label{sec:conclusions}

This work presents a digital twin simulation framework for pastillation processes that integrates a computer vision-based soft sensor (\texttt{PaNet}) into an automatic control system. The simulator produces realistic image data of pastilles on a conveyor belt, modeling flow rate and temperature fluctuations arising from random clogging events. Using this dataset, \texttt{PaNet-1D} (based on 1D convolution) and \texttt{PaNet-2D} (based on 2D convolution) were developed to predict the average temperature of the earliest row of pastilles, leveraging thermal image data for real-time monitoring. The conveyor belt speed was regulated through a PID control scheme to achieve precise temperature tracking and maintain production efficiency.
\\

\texttt{PaNet-1D} achieved high predictive accuracy, with RMSE values below 1.3~$^\circ$F for temperature prediction, and 0.4~$/\Delta t$ for flow rate prediction, outperforming \texttt{PaNet-2D}. Saliency map analysis showed that \texttt{PaNet-1D} effectively focused on the critical regions of the thermal images (the earliest row of pastilles), which underpinned its robust predictive capabilities. The closed-loop PID control system converged rapidly to a range of temperature setpoints, though minor discrepancies appeared between true and predicted temperatures during transient periods. In addition, fluctuations in the belt speed profile highlighted potential areas for improvement in the control strategy.
\\
 
Future work will explore advanced control algorithms, such as adaptive or model-predictive controllers, to mitigate speed fluctuations and further enhance system stability. Expanding the simulator to include more complex process dynamics and real-world noise will improve its applicability. Integration of \texttt{PaNet} with additional sensors, along with optimization of the control system for multi-objective performance metrics including energy efficiency and product quality, offers a promising direction for advancing automated pastillation technology. This framework illustrates the potential of computer vision and machine learning in driving progress in process control across diverse industrial applications.

\section{Acknowledgments} \label{sec:acknowledgements}

We acknowledge the support of the members of the Texas-Wisconsin-California Control Consortium and partial support from the National Science Foundation under grant CBET-2315963.

\bibliography{ref}
\bibliographystyle{unsrt}

\end{document}